\documentclass[10pt,a4paper,reqno]{amsart}

\usepackage[centertags]{amsmath}

\usepackage{latexsym,amssymb}
\usepackage[ansinew]{inputenc}
\usepackage[T1]{fontenc}
\usepackage{times,mathptmx}
\usepackage{color,geometry}
\usepackage{hyperref}

\usepackage{amsfonts, times, mathptmx}
\usepackage{amssymb}
\usepackage{amsthm}
\usepackage{newlfont}
\usepackage{url}
\usepackage{tikz}
\usepackage[english]{babel}

\newcommand {\real} {\mathbb{R}}

\theoremstyle{plain}
\newtheorem{thm}{Theorem} [section]

\newtheorem{lem}[thm]{Lemma}
\newtheorem{prop}[thm]{Proposition}
\newtheorem{conj}[thm]{Conjecture}

\theoremstyle{definition}

\newtheorem{rem}[thm]{Remark}

\begin {document}

\title{Reverse Cheeger inequality for planar convex sets}

\author{Enea Parini}

\date{\today}

\address{Aix Marseille Universit\'{e}, CNRS, Centrale Marseille, I2M, UMR 7373, 13453 Marseille, France}

\email{enea.parini@univ-amu.fr}



\thanks{}

\begin {abstract}
We prove the sharp inequality
\[ J(\Omega) := \frac{\lambda_1(\Omega)}{h_1(\Omega)^2} < \frac{\pi^2}{4},\]
where $\Omega$ is any planar, convex set, $\lambda_1(\Omega)$ is the first eigenvalue of the Laplacian under Dirichlet boundary conditions, and $h_1(\Omega)$ is the Cheeger constant of $\Omega$. The value on the right-hand side is optimal, and any sequence of convex sets with fixed volume and diameter tending to infinity is a maximizing sequence. Morever, we discuss the minimization of $J$ in the same class of subsets: we provide a lower bound which improves the generic bound given by Cheeger's inequality, we show the existence of a minimizer, and we give some optimality conditions.
\end{abstract}

\maketitle

\section{Introduction}

A celebrated inequality proven by Jeff Cheeger (\cite{cheeger}) states that, for every bounded domain $\Omega \subset \real^n$,
\begin{equation} \label{cheegerinequality} J(\Omega):=\frac{\lambda_1(\Omega)}{h_1(\Omega)^2} \geq \frac{1}{4}.\end{equation}
Here $\lambda_1(\Omega)$ is the first eigenvalue of the Laplacian under Dirichlet boundary conditions, while $h_1(\Omega)$ is the \emph{Cheeger constant} of $\Omega$, which is defined as
\[
h_1(\Omega)=\inf_{S\subset\Omega}\, \frac{P(S)}{|S|},
\]
$|S|$ being the $n$-dimensional Lebesgue measure of $S$, and $P(S)$ the distributional perimeter measured with respect to $\real^n$. A \emph{Cheeger set} is a set for which the infimum is attained. We refer to \cite{parinicheeger} for an introduction to the Cheeger problem.

The quantity $J$, up to an exponent, can be seen as a particular case of functionals of the kind
\[ 
\mathcal{R}_{p,q}(\Omega)=\frac{\lambda_1(p;\Omega)^{\frac{1}{p}}}{\lambda_1(q;\Omega)^{\frac{1}{q}}}, 
\]
where $\lambda_1(p;\Omega)$ is the first eigenvalue of the $p$-Laplacian under Dirichlet boundary conditions, with $1 \leq q < p \leq + \infty$. It is well known that $\lambda_1(p;\Omega) \to h_1(\Omega)$ as $p \to 1$ (see \cite{kawohlfridman}), while $\lambda_1(p;\Omega)^{\frac{1}{p}} \to \Lambda_1(\Omega)$ as $p \to +\infty$, where $\Lambda_1(\Omega)$ is the inverse of the radius of the biggest ball contained in $\Omega$ (see \cite{juutinenlindqvistmanfredi}). Some particular cases have already been considered in the literature. If $\Omega$ is a convex subset of $\real^2$, $p=\infty$ and $q=1$, it can be proved that
\begin{equation} \label{generalcase1} \frac{1}{2} \leq \frac{\Lambda_1(\Omega)}{h_1(\Omega)} < 1, \end{equation}
where the bounds are sharp. The first inequality becomes an equality when $\Omega$ is a ball, while a maximizing sequence is given by rectangles of the form $(-M,M)\times (-1/M,1/M)$. Let us now consider the case $p=\infty$ and $q=2$. Up to an exponent, the functional is equivalent to
\[ \frac{\Lambda_1(\Omega)^2}{\lambda_1(\Omega)}.\]
By the results in \cite[Section 8]{hersch} we have that, for any convex, planar domain $\Omega$,
\begin{equation} \label{generalcase2} \frac{1}{\lambda_1(B)} \leq \frac{\Lambda_1(\Omega)^2}{\lambda_1(\Omega)} < \frac{4}{\pi^2}.\end{equation}
The minimum is clearly attained when $\Omega$ is a ball, while the same sequence of elongating rectangles provides again a maximizing sequence.

In the spirit of inequalities \eqref{generalcase1} and \eqref{generalcase2}, we are interested in finding an upper bound for $J$. We are able to prove the inequality
\begin{equation} \label{reversecheeger} \frac{\lambda_1(\Omega)}{h_1(\Omega)^2} < \frac{\pi^2}{4}, \end{equation}
which can be seen as a \emph{reverse Cheeger inequality}. We mention that similar inequalities were considered by Buser (\cite{buser}) in the context of Riemannian geometry. We prove that inequality \eqref{reversecheeger} is sharp, and any sequence of planar, convex sets with fixed volume and diameter tending to infinity is a maximizing sequence. 

We then turn to the minimization of $J$ in the same class of sets. It can be observed that for typical choices of domains, such as circles or polygons, the value of $J$ is much higher than the lower bound given by \eqref{cheegerinequality}; therefore it is natural to wonder whether the inequality can be improved. Exploiting the knowledge of the behaviour of $J$ along sequences of sets with diverging diamenter, we are able to prove the existence of a minimizing set. However, simple computations show that the ball is no longer a minimizer, since a square or an equilateral triangle provide lower values of $J$; therefore, the identification of the optimal shape is no longer a trivial task. With this respect, we can provide partial results: the boundary of any minimizing set must be polygonal where it does not coincide with the boundary of its Cheeger set, in a sense which is explained in more detail in the following. It is still an open question to understand whether minimizing sets must actually be polygons.

The paper is structured as follows: after stating some preliminary results, we prove the continuity of the functional in the class of convex sets with respect to the Hausdorff distance (Section 3). In Section 4 we deal with the behaviour of $J$ along sequences of sets with diverging diameter and with the non-existence of a maximizer, thus proving inequality \eqref{reversecheeger}, while in Section 5 we prove the existence of a minimizer and we give some optimality conditions. Finally, we state some conjectures and open problems. An appendix about inequalities \eqref{generalcase1} and \eqref{generalcase2} complements the paper.

The author would like to thank Lorenzo Brasco for pointing out this problem to his attention and for many useful and interesting discussions.

\section{Preliminary results}
Let $\Omega \subset \real^n$ be an open set. The \emph{perimeter} of a set $E \subset \Omega$ (measured with respect to $\real^n$) is defined as
\[ P(E;\real^n):= |D\chi_E|(\real^n),\]
where $\chi_E$ is the characteristic function of $E$, and $|D\chi_E|(\real^n)$ is its total variation (see \cite{giusti}). The \emph{Cheeger constant} of $\Omega$ is
\[ h_1(\Omega) := \inf_{E \subset \Omega} \frac{P(E;\real^n)}{|E|},\]
where $|E|$ stands for the $n$-dimensional Lebesgue measure of $E$. A \emph{Cheeger set} is a set $C \subset \Omega$ such that
\[ \frac{P(C;\real^n)}{|C|} = h_1(\Omega). \]
The existence of a Cheeger set for every bounded Lipschitz domain $\Omega$ can be proved via the direct method of the Calculus of Variations. Uniqueness does not hold in general; however, any convex body has a unique Cheeger set, which is convex and with boundary of class $C^{1,1}$ (see \cite{altercaselles}). If $C$ is a Cheeger set for $\Omega$, then $\partial C \cap \partial \Omega$ is analytic, up to a closed singular set of Hausdorff dimension $n-8$; at the regular points of $\partial C \cap \Omega$, the mean curvature is equal to $\frac{h_1(\Omega)}{n-1}$ (see e.g. \cite[Proposition 4.2]{parinicheeger}).

If $\Omega$ is a convex, planar set, the Cheeger set can be characterized as the Minkowski sum of the so-called \emph{inner Cheeger set} and of a ball of radius $1/h_1(\Omega)$ (see \cite{kawohllachand}). In particular, $\partial C \cap \Omega$ consists of circular arcs.

It is interesting to observe that the Cheeger constant can be obtained as the limit for $p \to 1$ of the first eigenvalue of the $p$-Laplacian operator under Dirichlet boundary conditions, and it can be seen as the first eigenvalue of the $1$-Laplacian (see \cite{kawohlfridman}).

If $\Omega$ is a convex planar set, $C$ its Cheeger set, and $V \in C^1(\real^n;\real^n)$ a diffeomorphism, then the following shape derivative formula holds true (see \cite{parinisaintier}):
\begin{equation}\label{derh}
 h_1(\Omega;V)' = \frac{1}{|C|}\int_{\partial C \cap \partial \Omega} (\kappa -h_1(\Omega)) \langle V,\nu \rangle \,d\sigma, 
\end{equation}
where $h_1(\Omega;V)'$ is defined as
\[ h_1(\Omega;V)' := \lim_{t \to 0} \frac{h_1((Id+tV)(\Omega)) - h_1(\Omega)}{t},\]
$\kappa(x)$ is the curvature of $\partial \Omega$ at the point $x$, and $\nu$ is the unit exterior normal to $\partial\Omega$. For the reader's convenience, we also recall Hadamard's formula for the shape derivative of the first eigenvalue of the Laplacian $\lambda_1(\Omega)$:
\begin{equation} \label{derl} \lambda_1(\Omega;V)' = -\int_{\partial\Omega} |u_n|^2 \langle V,\nu \rangle\,d\sigma, \end{equation}
where $u_n$ is the normal derivative of $u$.

In the following we will need to know the explicit values of $\lambda_1(\Omega)$ and $h_1(\Omega)$ for some particular domains $\Omega$, which are listed in Table \ref{tabella1} (see \cite{siudeja} and \cite{kawohllachand}).

\begin{center}
\begin{table}[]
\begin{tabular}{|c|c|c|}  \hline $\Omega$ & $\lambda_1(\Omega)$ & $h_1(\Omega)$  \\ \hline triangle of area $A$ and perimeter $L$ & $\frac{\pi^2 L^2}{16A^2} < \lambda_1(\Omega) \leq \frac{\pi^2 L^2}{9A^2}$ & $\frac{L+\sqrt{4\pi A}}{2A}$ \\ \hline equilateral triangle & $\frac{\pi^2 L^2}{9A^2}$ & $\frac{L+\sqrt{4\pi A}}{2A}$ \\ \hline rectangle $(0,a)\times (0,b)$ & $\pi^2 \left(\frac{1}{a^2}+\frac{1}{b^2}\right)$ & $\frac{4-\pi}{a + b - \sqrt{(a - b)^2 + \pi ab}}$ \\ \hline circle of radius $R$  & $\frac{\lambda_1(B)}{R^2} \simeq \frac{5.7830}{R^2}$ & $\frac{2}{R}$ \\ \hline 
\end{tabular}
\vspace{4mm}
\caption{Explicit values of $\lambda_1(\Omega)$ and $h_1(\Omega)$ in some special cases.}
 \label{tabella1}
\end{table}
\end{center}

\section{Continuity of the functional}

Let $(\mathcal{K}_n,d_H)$ be the metric space given by the set of all open convex subsets of $\real^n$, endowed with the Hausdorff distance
\[ d_H(A,B):= d^H(\real^n\setminus A, \real^n \setminus B),\]
where
\[d^H(E,F):=\inf \{ \varepsilon \geq 0 \,|\,E \subset F_\varepsilon \text{ and }F \subset E_\varepsilon\} \]
and, for a set $E$,
\[ E_\varepsilon := \{ x \in \real^n \,|\,\text{dist}(x,E) \leq \varepsilon\}.\]
In this section we will prove that the functional
\[ J(\Omega):= \frac{\lambda_1(\Omega)}{h_1(\Omega)^2} \]
is continuous in $\mathcal{K}_n$ with respect to $d_H$. Since this fact is true for $\lambda_1(\Omega)$ (see for instance \cite[Theorem 2.3.17]{henrot}), we will prove that the claim holds for $h_1(\Omega)$. To prove our result, we will make use of the notion of $\Gamma$-convergence. Given a family of functionals $F_k : X \to \real$, where $X$ is a metric space, and given a limit functional $F_0 : X \to \real$, we say that $F_k \stackrel{\Gamma}{\to} F_0$ if the following conditions are satisfied:
\begin{itemize}
 \item \emph{liminf inequality:} for every sequence such that $x_k \to x$ in $X$, it holds
 \[ F_0(x) \leq \liminf_{k \to \infty} F_k (x_k).\]
 \item \emph{limsup inequality:} for every $x \in X$, there exists a sequence such that $x_k \to x$ in $X$ and
 \[ F_0(x) \geq \limsup_{k \to \infty} F_k (x_k)\]
 or, equivalently,
 \[ F_0(x) = \lim_{k \to \infty} F_k (x_k).\]
\end{itemize}
Suppose that $\{\overline{x}_k\}$ is such that $F_k(\overline{x}_k)=\min_X F_k$. It can be proved that, if there exists a compact set $E \subset X$ such that $\overline{x}_k \in E$ for $k$ sufficiently big, then $\overline{x}_k \to \overline{x}$ (up to a subsequence), where $F(\overline{x})=\min_X F$, and $\min_X F_k \to \min_X F$ (see \cite{braides}).

\begin {prop} \label{continuityh1}
Let $\Omega,\,\Omega_k\subset\real^n$ be bounded open convex sets such that $\Omega_k \to \Omega$ in the Hausdorff metric. Then,
\[
\lim_{k\to\infty} h_1(\Omega_k)=h_1(\Omega).
\]
and the corresponding Cheeger sets converge in the Hausdorff metric.
\end{prop}

\begin{proof}
From the boundedness of $\Omega$, the convexity of the sequence $\{\Omega_k\}_{k\in\mathbb{N}}$ and the convergence, we get that there exists a (convex) set $F$ such that
\[
\Omega\subset F \qquad \mbox{ and }\qquad \Omega_k\subset F,\quad k\in\mathbb{N}.
\]
Let
$\mathcal{K}_n(\Omega_k)$, $\mathcal{K}_n(\Omega)$ and $\mathcal{K}_n(F)$ be the families of open convex subsets of
$\Omega_k$, $\Omega$ and $F$ respectively. Let us define the functionals
\[\Phi_k(E):=\left\{\begin{array}{c l}\displaystyle\frac{P(E)}{|E|}&\textrm{for }E\in\mathcal{K}_n(\Omega_k) \\ +\infty & \text{for }E \in \mathcal{K}_n(F) \setminus \mathcal{K}_n(\Omega_k)\end{array}\right.\] and
\[\Phi(E):=\left\{\begin{array}{c l}\displaystyle\frac{P(E)}{|E|}&\textrm{for }E\in\mathcal{K}_n(\Omega) \\ +\infty & \text{for }E \in \mathcal{K}_n(F) \setminus \mathcal{K}_n(\Omega)\end{array}\right.\] Observe that
\[h_1(\Omega_k)=\inf_{E\in\mathcal{K}_n(\Omega_k)} \frac{P(E)}{|E|} = \inf_{E\in\mathcal{K}_n(F)} \Phi(E)\]
since every convex domain admits a unique convex Cheeger set. We are now ready to prove
the $\Gamma$-convergence of the functionals $\Phi_k$ to $\Phi$.\\
\textit{liminf inequality}. Let $E, E_k \in \mathcal{K}_n(F)$ be
such that $E_k\to E$ in the Hausdorff metric. This implies that $E_k \to E$ in the $L^1$-topology. If $E_k \in \mathcal{K}_n(\Omega_k)$ only for a finite number of elements, then there is nothing to prove. Otherwise, we have that $E \in \mathcal{K}_n(\Omega)$. Of course we have
$|E_k|\to |E|$, while from the lower semicontinuity of the
perimeter  we obtain
$\displaystyle P(E)\leq \liminf_{k\to\infty} P(E_k)$. In
conclusion we get
\[\Phi(E)\leq\liminf_{k\to\infty} \Phi_k(E_k).\]\\
\textit{limsup inequality}. Let $E\in\mathcal{K}_n(F)$; if $E \not\in \mathcal{K}_n(\Omega)$ there is nothing to prove. Let us now suppose $E \in \mathcal{K}_n(\Omega)$, and let us define
$E_k:=E\cap \Omega_k$. The sets $E_k$ are convex sets contained in
$\Omega_k$, and are such that $E_k \to E$ in the Hausdorff metric (see \cite[p. 32]{henrotpierre}), and therefore in the $L^1$-topology.
From \cite[Lemma 4.4]{buttazzoferone} one has $P(E_k)\to P(E)$,
so that
\[\Phi(E)=\lim_{k\to\infty} \Phi_k(E_k).\]
\textit{Equicoercivity}. Let $C_k$ be a convex Cheeger
set for $\Omega_k$. The sets $C_k$ are all contained in $F$, and therefore they are elements of $\mathcal{K}_n(F)$, which is a compact set.\\
From the properties of the $\Gamma$-convergence we obtain that,
after possibly passing to a subsequence,
\[h_1(\Omega_k)\to h_1(\Omega)\]
and there exists a sequence of Cheeger sets $C_k$ for
$\Omega_k$ converging in the Hausdorff metric to a Cheeger set
$C$ for $\Omega$. However, by uniqueness of the limit we have that the whole sequences $h_1(\Omega_k)$ and $C_k$ are convergent.
\end{proof}

\begin{prop} \label{continuityJ}
The functional $J$ is continuous in $\mathcal{K}_n$ with respect to $d_H$. 
\end{prop}
\begin{proof}
The claim follows from Proposition \ref{continuityh1} and from the continuity of $\lambda_1(\Omega)$ (see for instance \cite[Theorem 2.3.17]{henrot}).
\end{proof}

\section{Behaviour at infinity and non-existence of a maximizer}

In the following we will restrict ourselves to the two-dimensional case. For the sake of simplicity, we set $\mathcal{K}:=\mathcal{K}_2$. In this section we investigate the behaviour of $J(\Omega)$ along sequences $\{\Omega_k\}$ of convex planar sets such that $|\Omega_k|$ is fixed and $\text{diam }\Omega_k \to +\infty$. This result will be crucial in the next section in order to prove existence of a minimizer in $\mathcal{K}$. Moreover, we will see that every such sequence is a maximizing sequence, and that the supremum is not attained. In the proofs we will use the basic, but important observation that the functional $J$ is invariant by rigid motions and by dilations, since both $\lambda_1(\Omega)$ and $h_1(\Omega)^2$ have the same scaling: if we define for $t>0$
\[ t\,\Omega:=\{ x \in \real^2 \,|\,t^{-1}\,x \in \Omega \},\]
then
\[ \lambda_1(t\,\Omega)= t^{-2}\,\lambda_1(\Omega) \qquad \text{and}\qquad h_1(t\,\Omega)=t^{-1}\,h_1(\Omega).\]
If $\Omega \subset \real^2$, we can define its \emph{Steiner symmetrization} $\Omega^*$ with respect to the $x$-axis. Let the symbol $\mathcal{H}^1$ stand for the one-dimensional Hausdorff measure. For $t \in \real$, define $m(t):=\mathcal{H}^1(\Omega \cap \{x=t\})$, and let $I(x)=\left( -\frac{m(x)}{2},\frac{m(x)}{2} \right)$ if $m(x)>0$, or $I(x)=\emptyset$ if $m(x)=0$. Then,
\[ \Omega^* :=  \bigcup_{x \in \real} \left(\{x\} \times I(x)\right).\]
It is well-known that the first eigenvalue of the Laplacian under Dirichlet boundary conditions decreases under symmetrizations, so that $\lambda_1(\Omega^*) \leq \lambda_1(\Omega)$ (see for instance \cite[Theorem 2.2.4]{henrot}).

\begin{prop} \label{behaviourinfinity}
Let $\Omega \in \mathcal{K}$. Then,
\[ \frac{\lambda_1(\Omega)}{h_1(\Omega)^2} < \frac{\pi^2}{4}.\]
Moreover, every sequence $\{\Omega_k\}$ in $\mathcal{K}$ such that $|\Omega_k|=V$ for some $V>0$, and $\text{diam }\Omega_k \to + \infty$ as $k \to \infty$ satisfies
\[ \frac{\lambda_1(\Omega_k)}{h_1(\Omega_k)^2} \to \frac{\pi^2}{4} \]
as $k \to \infty$.
\end{prop}
\begin{proof}
Let $\Omega \in \mathcal{K}$, and let $C$ be its Cheeger set, which is a convex set. By monotonicity, $\lambda_1(\Omega)\leq \lambda_1(C)$, while $h_1(\Omega) = \frac{|\partial C|}{|C|}$. By the results in \cite{polya} it holds
\begin{equation} \label{sup} \frac{\lambda_1(\Omega)}{h_1(\Omega)^2} \leq \frac{\lambda_1(C)|C|^2}{|\partial C|^2} < \frac{\pi^2}{4},\end{equation}
Now we prove the second part of the claim. Since the functional $J$ is invariant by rotations, it is possible to rotate the sets $\Omega_k$ in such a way that the set $(0,d_k)\times \{0\}$ is contained in $\Omega_k$. Define $\varepsilon_k$ as the ``inner width'' of $\Omega_k$, that is,
\[ \varepsilon_k = \max \left\{\varepsilon > 0 \,\bigg|\,\varepsilon=b-a,\, \{x\}\times (a,b) \subset \Omega_k \text{ for some }x \in (0,d_k)\right\}.\]
It is then clear that $\varepsilon_k \leq \frac{2V}{d_k}$, because otherwise $\Omega_k$ would contain a quadrilateral of volume bigger than $V$, and therefore $\varepsilon_k \to 0$ as $k \to \infty$. $\Omega_k$ contains a triangle $T_k$ of basis $\varepsilon_k$ and height $a_k d_k$, with $a_k \in \left[\frac{1}{2},1\right]$. Let $\Omega_k^*$ be the Steiner symmetrization of $\Omega_k$ with respect to the $x$-axis. $\Omega_k^*$ is also convex, and it is contained in a rectangle $R_k$ of edges $d_k$ and $\varepsilon_k$. After passing to a subsequence, one can suppose that $a_k \to a$. But then (see Table \ref{tabella1})
\[ \lambda_1(\Omega_k) \geq \lambda_1(\Omega_k^*) \geq \lambda_1(R_k) = \pi^2 \left(\frac{1}{d_k^2}+ \frac{1}{\varepsilon_k^2}\right)\]
\[ h_1(\Omega_k) \leq h_1(T_k) \leq \frac{\varepsilon_k + a_kd_k + \sqrt{\varepsilon_k^2 + (a_kd_k)^2} + \sqrt{2\pi \varepsilon_k a_k d_k}}{\varepsilon_k a_k d_k}\]
and therefore
\begin{equation} \label{liminf} \frac{\lambda_1(\Omega_k)}{h_1(\Omega_k)^2}\geq \frac{\lambda_1(R_k)}{h_1(T_k)^2} \Rightarrow \liminf_{k \to \infty}\frac{\lambda_1(\Omega_k)}{h_1(\Omega_k)^2} \geq \frac{\pi^2}{4}.\end{equation}
By \eqref{liminf} and \eqref{sup} we obtain the claim.
\end{proof}

\begin{rem}
Observe that for every interval $I=[a,b]$, there holds $\lambda_1(I)=\frac{\pi^2}{(b-a)^2}$ and $h_1(I)=\frac{2}{b-a}$, so that
\[ \frac{\lambda_1(I)}{h_1(I)^2}=\frac{\pi^2}{4}.\]
\end{rem}

\section{Existence of a minimizer and optimality conditions}

In this section we will prove the existence of a minimizer for $J$ in the class $\mathcal{K}$ by means of the direct method of the Calculus of Variations; to this end, an essential tool will be Proposition \ref{behaviourinfinity}. A key observation is that the ball does not minimize $J$, but it is a critical point. Moreover, we will derive some optimality conditions.

It is useful to compute some explicit values of $J$. In Table \ref{tabella2} we give the values of $J(\Omega)$ where $\Omega$ is a regular $n$-gon (with edge length equal to $1$) or the unit circle. The values are computed analytically (where possible) or numerically. For the value of $\lambda_1(\Omega)$ in a hexagon see \cite{curetonkuttler}, while a formula for $h_1(\Omega)$ when $\Omega$ is a regular polygon was given in \cite{kawohllachand}.

\begin{center}
\begin{table}[]
\begin{tabular}{|c|c|c|c|}  \hline $n$ & $\lambda_1(\Omega)$ & $h_1(\Omega)$ & $J(\Omega)$ \\ \hline 3 & 52.63789 & 6.157649 & 1.388252 \\ \hline 4 & 19.739208 & 3.772453 & 1.38701 \\ \hline 5 & 10.9964  & 2.8044 & 1.39820 \\ \hline 6 & 7.15533 & 2.2543 & 1.40801 \\ \hline 8 & 3.7988 & 1.6351 & 1.42088  \\  \hline $\infty$ & 5.7830 & 2 & 1.4457 \\ \hline 
\end{tabular}
\vspace{4mm}
\caption{Explicit values of $\lambda_1(\Omega)$, $h_1(\Omega)$ and $J(\Omega)$ for some regular $n$-gons (edge length = 1) and for the unit circle ($n=\infty$).}
\label{tabella2}
\end{table}
\end{center}

\medskip

\begin{prop} \label{nonoptimallowerbound}
For every $\Omega \in \mathcal{K}$ it holds
\[ \frac{\lambda_1(\Omega)}{h_1(\Omega)^2} \geq \frac{\pi^2}{16}\, (\simeq 0.616...). \]
\end{prop}
\begin{proof}
The claim follows from the fact that
\[ \frac{\lambda_1(\Omega)}{h_1(\Omega)^2}=\frac{\Lambda_1(\Omega)^2}{h_1(\Omega)^2}\cdot \frac{\lambda_1(\Omega)}{\Lambda_1(\Omega)^2},\]
and from the results in \cite[Section 8]{hersch} and in Appendix \ref{appcheeger}. 
\end{proof}

We remark that the lower bound in Proposition \ref{nonoptimallowerbound} is non-optimal, since otherwise $\Omega$ should be at the same time a ball and an infinite strip, but it improves the generic lower bound $\frac{1}{4}$ given by Cheeger's inequality.

\begin{prop} \label{existence}
The functional $J$ admits a minimizer in the class $\mathcal{K}$ of planar convex sets. 
\end{prop}
\begin{proof}
Let $\{\Omega_k\}$ be a minimizing sequence. Since the functional is scaling invariant, without loss of generality we can suppose that all the sets $\Omega_k$ have the same volume $V$. We will prove that $d_k:=\text{diam } \Omega_k$ is uniformly bounded. Suppose by contradiction that this is not the case; hence, we can suppose that there exists a subsequence (still denoted by $\Omega_k$) such that $d_k=\text{diam } \Omega_k \to \infty$ as $k\to \infty$. By Proposition \ref{behaviourinfinity},
\[ \lim_{k \to \infty} \frac{\lambda_1(\Omega_k)}{h_1(\Omega_k)^2} = \frac{\pi^2}{4}, \]
a contradiction to the fact that $\Omega_k$ is a minimizing sequence, since for a ball $B$
\[
\frac{\lambda_1(B)}{h_1(B)^2} \simeq 1.4457 < \frac{\pi^2}{4}. 
\] As a consequence, $d_k:=\text{diam } \Omega_k$ is uniformly bounded. Then, up to a suitable translation, there exists a ball containing all the sets $\Omega_k$. Therefore, one can extract a convergent subsequence, and the claim follows from Proposition \ref{continuityJ}. 
\end{proof}

In the following we will obtain some optimality conditions for a minimizer $\Omega$. If $C$ is its Cheeger set, we will analyze separately the part of $\partial \Omega$ ``far'' from the Cheeger set, that is $\partial \Omega \setminus \partial C$, and $\partial \Omega \cap \partial C$, which is of class $C^{1,1}$.

\begin{prop} \label{notaball}
The ball does not minimize $J$. 
\end{prop}
\begin{proof}
This follows from the easy observation that, for a square $Q$ and a ball $B$, $J(Q)<J(B)$. 
\end{proof}

\begin{lem} \label{lemmasamecheeger}
Let $\Omega \subset \real^2$ be a convex set, and let $C$ be its Cheeger set. Let $\Omega' \subset \real^2$ be a convex set such that $\Omega \subset \Omega'$ and $\partial \Omega' \cap \partial C = \partial \Omega \cap \partial C$. Then $h_1(\Omega')=h_1(\Omega)$, and $C$ is the Cheeger set of $\Omega'$. 
\end{lem}
\begin{proof}
It is known (see \cite[Theorem 1]{kawohllachand}) that $C$ is given by the Minkowski sum $C_r \oplus B_r$, where $r=h_1(\Omega)^{-1}$, $B_r$ is the ball of radius $r$, and $C_r$ is the \emph{inner Cheeger set} defined as
\[ C_r := \{ x \in \Omega\,|\,\text{dist}(x;\partial \Omega) \geq r\}.\]
Moreover, $r$ is the unique value such that $|C_r|=\pi r^2$. We set $\Gamma_1 := \partial \Omega \cap \partial C$, $\Gamma_2 = \partial \Omega \setminus \partial C$ and $\Gamma_3 = \partial \Omega' \setminus \partial C$. Define
%$\Gamma_1 := \{\}$, $\Gamma_2 = \partial \Omega \setminus \partial C$ and $\Gamma_3 = \partial \Omega' \setminus \partial C$. Define
\[ C'_r := \{ x \in \Omega\,|\,\text{dist}(x;\partial \Omega') \geq r\}.\]
We want to prove that $C'_r = C_r$. We have $C_r \subset C'_r$. Moreover, for every $x \in \partial C_r$, we have that $\text{dist}(x;\partial \Omega)=r$, which means that there exists $y \in \partial \Omega$ such that $\text{dist}(x,y)=r$; since every point in $\partial B_r(x)$ belongs to $\overline{C}$, it holds $y \in \partial C$, and therefore $\text{dist}(x;\Gamma_1)=r$. Since $\text{dist}(x;\Gamma_3)\geq \text{dist}(x;\Gamma_2) \geq r$, one has $\text{dist}(x;\partial \Omega')=r$, which means that $x \in \partial C_r'$. Since $C_r$ and $C_r'$ are convex sets, we must have $C_r=C_r'$. From the definition of $C'_r$ and the fact that $|C'_r|=\pi r^2$, it follows that $C$ is the Cheeger set of $\Omega'$, and $h_1(\Omega')=h_1(\Omega)$.
\end{proof}

\begin{center}
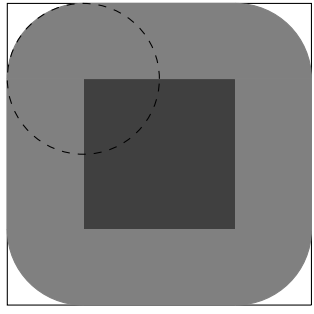
\begin{figure}
\begin{tikzpicture}
\draw (4,0)--(4,4)--(0,4)--(0,0)--(4,0);
\draw[dashed] (3,1)--(3,3)--(1,3)--(1,1)--(3,1);
\filldraw[fill=darkgray] (1,1) rectangle(3,3);
\draw[fill=gray, draw=gray] (0,3)--(1,3)--(1,4) arc(90:180:1cm);
\draw[fill=gray, draw=gray] (3,4)--(3,3)--(4,3) arc(0:90:1cm);
\draw[fill=gray, draw=gray] (1,0)--(1,1)--(0,1) arc(180:270:1cm);
\draw[fill=gray, draw=gray] (4,1)--(3,1)--(3,0) arc(270:360:1cm);
\draw[fill=gray, draw=gray] (1,0) rectangle(3,1);
\draw[fill=gray, draw=gray] (0,1) rectangle(1,3);
\draw[fill=gray, draw=gray] (1,3) rectangle(3,4);
\draw[fill=gray, draw=gray] (3,1) rectangle(4,3);
\draw[dashed, draw=black] (2,3) arc(0:360:1cm);

\end{tikzpicture}
\caption{The Cheeger set $C$ of a convex set $\Omega$ (gray) is given by the Minkowski sum of the inner Cheeger set (dark gray) and of a ball of radius $r=h_1(\Omega)^{-1}$.}
\end{figure}
\end{center}

\begin{prop} \label{fuoridacheeger}
Let $\Omega$ be a minimizer. Let $C$ be the Cheeger set of $\Omega$, and let $F:=\partial \Omega \setminus \partial C$. Then each connected component of $F$ consists of two segments intersecting at a common vertex. 
\end{prop}
\begin{proof}
Let $\Gamma$ be a connected component of $F$. The closure of $\Gamma$ in the relative topology intersects $\partial C$ in two distinct points $A$ and $B$ where $\Omega$ admits a tangent line: this is due to the convexity of $\Omega$ and to the fact that, thanks to the characterization of Cheeger sets, there exists an interior tangent ball at both $A$ and $B$. Let $t_A$ and $t_B$ the tangents at $A$ and $B$ respectively, and suppose that they intersect at a point $P$. Define $\Omega'$ as the interior of the convex hull of $\{P\} \cup \Omega$, and suppose by contradiction that $\Omega \neq \Omega'$. By monotonicity, $\lambda_1(\Omega')<\lambda_1(\Omega)$, while by Lemma \ref{lemmasamecheeger} one would have $h_1(\Omega)=h_1(\Omega')$, a contradiction to the fact that $\Omega$ minimizes $J$. The case where $t_A$ and $t_B$ are parallel can be ruled out in a similar way, by considering an arbitrary convex domain $\Omega' \supset \Omega$ such that $\partial \Omega$ and $\partial \Omega'$ coincide on $\partial \Omega \setminus \Gamma$. 
\end{proof}

Now we will obtain necessary conditions on $\partial \Omega \cap \partial C$ by using shape optimization arguments. If $V \in C_c^1(\real^n;\real^n)$, the shape derivative of $J$ in direction $V$, defined as
\[ J(\Omega;V)' := \lim_{t \to 0} \frac{J((Id+tV)(\Omega)) - J(\Omega)}{t} \]
exists due to the fact that both $\lambda_1(\Omega)$ and $h_1(\Omega)$ are shape differentiable (see \cite{parinisaintier}). We observe that $\Omega$ is a critical point of the functional $J$ if and only if
\begin{equation} \label{shapederivative} \lambda_1(\Omega;V)' = \frac{2\lambda_1(\Omega)}{h_1(\Omega)} h_1(\Omega;V)'. \end{equation}
for every deformation $V$ preserving the convexity of $\Omega$. More precisely, by \cite[Theorem 7]{henrotoudet} the shape derivative $\lambda_1(\Omega;V)'$ can be written for every $V$ which modifies only the strictly convex parts of $\partial \Omega$. Arguing as in \cite[Theorem 7]{henrotoudet} and adapting \cite[Lemma 2.1]{coxross} for $h_1(\Omega)$, we see that the same holds true for $h_1(\Omega;V)'$. By the arbitrariness of $V$, and by relations \eqref{derh} and \eqref{derl}, this implies that on the strictly convex parts of the set $\partial \Omega \cap \partial C$
\begin{equation} \label{relation} |u_n|^2 = \frac{2 \lambda_1(\Omega)}{h_1(\Omega)\,|C|} (h_1(\Omega) - \kappa ) .\end{equation}

\begin{prop} \label{ballnotminimal}
The ball is a critical point of $J$.
\end{prop}
\begin{proof}
The following representation formula for $\lambda_1(\Omega)$ holds true (see \cite{rellich}):
\[ \lambda_1(\Omega) = \frac{1}{2}\int_{\partial \Omega} |u_n|^2\,\langle x, \nu \rangle\,d\mathcal{H}^{n-1},\]
where $u$ is a normalized eigenfunction, and $\nu$ is the outer normal vector. If $\Omega=B$ is the unit ball, then $|u_n|=c$ on $\partial B$, so that 
\[ c = \left(\frac{\lambda_1(B)}{\pi}\right)^{\frac{1}{2}}. \]
Since
\[ \frac{\lambda_1(B)}{\pi} = \frac{2 \lambda_1(B)}{2 \pi}(h_1(B) - \kappa(x))\]
for every $x \in \partial B$, we have that \eqref{relation} holds true and therefore that the ball is a critical point. However, by Proposition \ref{notaball} it is not a minimizer.
\end{proof}

Condition \eqref{relation} implies that a normalized first eigenfunction $u$ satisfies the partially overdetermined problem
\begin{equation} \label{overdetermined} \left\{\begin{array}{r c l l} -\Delta u & = & \lambda_1(\Omega)\,u & \text{in }\Omega \\ u & = & 0 & \text{on }\partial \Omega \\ |u_n|^2 & = & a-b\kappa & \text{on }\Gamma
   \end{array}\right. \end{equation}
where $a= \frac{2 \lambda_1(\Omega)}{|C|}$ and $b=- \frac{2 \lambda_1(\Omega)}{h_1(\Omega)\,|C|}$, and $\Gamma$ is the union of the strictly convex parts of $\partial \Omega \cap \partial C$; each component of $\Gamma$ is therefore of class $C^{1,1}$. However, it can be shown that $\Gamma$ actually enjoys higher regularity, as stated in the following result.
\begin{prop}
Each component of $\Gamma$ is of class $C^\infty$.
\end{prop}
\begin{proof}
Let $\Gamma'$ be a component of $\Gamma$. We already know that $\Gamma' \in C^{1,1}$. By elliptic regularity, $u \in C^{1,\alpha}(\Omega \cup \Gamma')$ (see \cite[Corollary 8.36]{gilbargtrudinger}). This implies that $|u_n| = |\nabla u| \in C^\alpha(\Gamma')$, and therefore $|u_n|^2 \in C^\alpha(\Gamma')$, which implies in turn that $\kappa \in C^\alpha(\Gamma')$. If $v : I \subset \real \to \real$ is the function whose graph describes $\Gamma'$ locally, it holds 
\[  - \left(\frac{v'}{\sqrt{1 + (v')^2}}\right)' = \kappa \]
in the weak sense. Since $v \in C^{1,1}$, $v'$ is bounded, and hence Schauder regularity results apply. This implies that $\Gamma' \in C^{2,\alpha}$, and $u \in C^{2,\alpha}(\Omega \cup \Gamma')$. By a bootstrap argument, we obtain that $\Gamma' \in C^\infty$, and $u \in C^\infty(\Omega \cup \Gamma')$. 
\end{proof}

Finally, we prove that $\partial \Omega \cap \partial C$ can not contain any arc of circle. This is a consequence of a recent result by Fragal\`{a} and Gazzola about partially overdetermined boundary value problems \cite{fragalagazzola}.

\begin{prop}
Let $\Omega$ be a minimizer. Then $\partial \Omega \cap \partial C$ can not contain arcs of circle. 
\end{prop}
\begin{proof}
Suppose by contradiction that $\partial \Omega \cap \partial C$ contains an arc of circle $\Gamma$. If $u$ is a normalized eigenfunction on $\Omega$, we have by \eqref{overdetermined} that $u_n$ is constant on $\Gamma$. By \cite[Theorem 1]{fragalagazzola}, $\Omega$ must be a ball, a contradiction to Proposition \ref{ballnotminimal}.
\end{proof}

\section{Final remarks and open problems}

Many questions concerning the maximization and the minimization of the functional $J$ remain open. For instance, it would be interesting to generalize the results to higher dimensions, and to prove a reverse Cheeger inequality also for non-convex sets. A major difficulty in these cases is the fact that a lot of information on Cheeger sets, such as explicit values for particular domains, uniqueness and regularity, is lacking.

As for the minimization of $J$ among planar, convex sets, it remains open to prove that any minimizer $\Omega$ is a polygon. Moreover, if $\Omega$ was a polygon, we do not know whether it should be regular; in this case, explicit computations (see Table \ref{tabella2}) support the claim that a minimizer should be a square. In any case, it is easy to show that the square minimizes $J$ among all rectangles. We mention a recent result by Bucur and Fragal\`{a}, which states that among all polygons of $n$ sides with fixed volume, $h_1(\Omega)$ is minimized for the regular one (see \cite{bucurfragala}). This fact is known to be true also for $\lambda_1(\Omega)$, if one restricts to the classes of triangles or quadrilaterals (see \cite{polyaszego}). However, in view of \cite{bucurfragala}, the claim that among all convex polygons with fixed number of edges, $J$ is minimized by the regular one, is actually stronger, at least in the subclass of convex sets, than the well-known P\'{o}lya-Szeg\H{o} conjecture, which 
states the same claim for $\lambda_1(\Omega)
$ and is still open for $n \geq 5$. 

Finally, one could wonder whether this kind of results holds true also for the more general functional $\mathcal{R}_{p,q}$. We state the following conjecture.

\begin{conj}
Functionals of the kind
\[ \mathcal{R}_{p,q}(\Omega) = \frac{\lambda_1(p;\Omega)^{\frac{1}{p}}}{\lambda_1(q;\Omega)^{\frac{1}{q}}} \]
with $q < p$ admit a minimizer, and they are bounded from above. A maximizing sequence is given by any sequence of sets of fixed volume such that $\text{diam }\Omega_k \to \infty$, the supremum is not attained and is equal to
\[ \frac{\lambda_1(p;I)^{\frac{1}{p}}}{\lambda_1(q;I)^{\frac{1}{q}}},\]
where $I \subset \real$ is an interval.
\end{conj}

The value of $\lambda_1(p;I)$ can be determined explicitly. If $I=(a,b)$ and $p \in (1,+\infty)$, then
\[ \lambda_1(p;I)= (p-1)\left( \frac{2\pi}{p(b-a)\sin{\left(\frac{\pi}{p}\right)}} \right)^p\]
(see \cite{otani}), so that
\[ \frac{\lambda_1(p;I)^{\frac{1}{p}}}{\lambda_1(q;I)^{\frac{1}{q}}} = \frac{q(p-1)^{\frac{1}{p}}}{p(q-1)^{\frac{1}{q}}}\cdot\frac{\sin{\left(\frac{\pi}{q}\right)}}{\sin{\left(\frac{\pi}{p}\right)}}.\]

\appendix

\section{Optimal inequalities for $\mathcal{R}_{\infty,1}$ and $\mathcal{R}_{\infty,2}$} \label{appcheeger}

In this section we will discuss the inequalities given in \eqref{generalcase1} and \eqref{generalcase2}, showing that the bounds are sharp. Let us first prove \eqref{generalcase1}. Recall that
\[ \mathcal{R}_{\infty,1}(\Omega)=\frac{\Lambda_1(\Omega)}{h_1(\Omega)}, \]
where $\Lambda_1(\Omega)$ is the first eigenvalue of the infinity Laplacian, that is, the inverse of the radius $R$ of the biggest ball contained in $\Omega$. It holds
\[ h_1(\Omega) \leq \frac{2}{R} = 2 \Lambda_1(\Omega) \Rightarrow \frac{\Lambda_1(\Omega)}{h_1(\Omega)} \geq \frac{1}{2},\]
and the infimum is attained for a ball. On the other hand, the radius of the biggest ball contained in $\Omega$ is bigger than the radius of the balls whose union is the Cheeger set for $\Omega$. This implies
\[ \Lambda_1(\Omega) \leq h_1(\Omega).\]
In fact, by the characterization of Cheeger sets in \cite[Theorem 1]{kawohllachand}, one has 
\[ \Lambda_1(\Omega) < h_1(\Omega). \]
Indeed, $C = C_r \oplus B_r$, where $r=\frac{1}{h_1(\Omega)}$, and $C_r$ is the ``inner Cheeger set'', which satisfies $|C_r|=\pi r^2$. If $r=R$, the radius of the biggest ball contained in $\Omega$, we would have that $C=B_R$ and hence $|C_r|=0$, a contradiction.
A sequence of rectangles of the form $R_d:=(-d,d)\times (-1,1)$ provides a maximizing sequence, since
\[ \frac{\Lambda_1(R_d)}{h_1(R_d)} \to 1 \]
as $d \to \infty$. Let us now come to \eqref{generalcase2}. In \cite[Section 8]{hersch} the inequality is stated as
\[ \frac{1}{\lambda_1(B)} \leq \frac{\Lambda_1(\Omega)^2}{\lambda_1(\Omega)} \leq \frac{4}{\pi^2}, \]
hence without strict inequality for the upper bound. But looking carefully at the proof, one can notice that in the case $\beta = \pi$ the first eigenvalue can be estimated strictly from below by the eigenvalue in the infinite strip, while in the case $\beta < \pi$, for $\lambda_1(T)$ the estimate
\[ \lambda_1(T) >  \frac{\pi^2}{4} \Lambda_1(\Omega)^2 \]
holds true (see \cite[Theorem 1.2]{siudeja}).

\bibliographystyle{amsplain}
\bibliography{biblioreverse}

\end {document}